\documentclass[11pt]{article}
\usepackage{latexsym,amsthm,verbatim,ifthen,graphicx,color,mathrsfs}
\usepackage[all]{xy}
\usepackage{color}
\usepackage{tikz}
\usepackage{xcolor}
\usepackage{tikz-cd}

\newcommand\blfootnote[1]{%
  \begingroup
  \renewcommand\thefootnote{}\footnote{#1}%
  \addtocounter{footnote}{-1}%
  \endgroup
}

\usepackage{color}

\usepackage{lscape}
 \usepackage[ansinew]{inputenc}
 \usepackage{multicol}
 \usepackage[all]{xy}\xyoption{poly}
 \usepackage{mathrsfs}
\usepackage[psamsfonts]{eucal}
\usepackage{amssymb, amsmath}
\usepackage{geometry}
\usepackage{enumerate}
\usepackage{float}

\usepackage{hyperref}  
\hypersetup{
    colorlinks=true,
    linkcolor=blue,
    filecolor=magenta,      
    urlcolor=cyan,
    citecolor=blue,
}
\usepackage{eucal}

\newtheorem*{theorem*}{Theorem}
\newtheorem{theorem}{Theorem}[section]

 \newtheorem{proposition}[theorem]{Proposition}
 \theoremstyle{definition}
 \newtheorem{definition}[theorem]{Definition}
 \theoremstyle{remark}
 \newtheorem{remark}[theorem]{Remark}
 \newtheorem{example}[theorem]{Example}

\newcommand{\im}{\operatorname{Im}}
\newcommand{\id}{\operatorname{id}}

\newcommand{\aut}{\operatorname{aut}}

\newcommand{\E}{{\mathcal E}}

\newcommand{\F}{\mathcal F}
\renewcommand{\H}{{\CMcal{H}}}

\newcommand{\ev}{{\text{ev}}}
\newcommand{\Fib}{\operatorname{\mathcal{F}\text{\normalfont ib}}}
\newcommand{\fib}{\Fib}
\newcommand{\FF}{\CMcal{F}}

\newcommand{\quie}{\E}

\begin{document}

\title{Classification of fiber sequences with a prescribed holonomy action}

\author{Mario Fuentes\footnote{The author has been partially supported by the MICINN grant MTM2016-78647-P and the grant FPU18/04140.}}

\maketitle

\blfootnote{
2020 {\em Mathematics Subject Classification}. 55R15, 55R05.}
\blfootnote{ {\em Key words and phrases}. Fiber spaces, classifying spaces, action of the fundamental group. }

\begin{abstract}
We define  $\H$-fibration sequences as fibrations where the holonomy action of the fundamental group of the base on the fiber lies in a given subgroup $\H$ of $\E(F)$, where $\E(F)$ is the homotopy automorphism group of the fiber. Furthermore, we classify these $\H$-fibration sequences via a universal $\H$-fibration sequence.
\end{abstract}

\section*{Introduction}

From the foundational paper of Stasheff \cite{sta} to the recent reference \cite{blomcha} of homotopy theoretical flavour, the literature is splashed of results which let us classify fibrations, sometimes of  a certain type (see for instance \cite{gargarmu,wa}), by means of a suitable  classifying space. The most general result in the
 homotopy category is probably the classical work of May \cite{may} which classifies fibrations with a given {\em category of fibers}. 
 
 In these works, it is constructed a classifying object $B\aut F$, and a universal quasi-fibration sequence
 $$F\to B(\ast,\aut F,F)\to B\aut(F)$$
 such that, given a space $B$ of the homotopy type of a CW-complex, there is a bijection
 $$\Lambda:[B,B\aut F]\to \Fib(B,F)$$
 where $[B,B\aut F]$ is the set of homotopy classes of maps from $B$ to $B\aut F$ and $\Fib(B,F)$ is the set of equivalence classes of fibration sequence over $B$ with fiber $F$.
 
 However, this treatment does not readily apply to classify fibrations with a prescribed holonomy action of the fundamental group of the base on the fiber: given a fibration sequence $F\stackrel{\omega}{\longrightarrow} E\stackrel{\pi}{\longrightarrow} B$, there is map of groups $\pi_1(B,b_0)\to \E(F)=\pi_0(\aut F)$. Given a subgroup $\H$ of $\E(F)$, we say that a fibration sequence is an $\H$-fibration sequence if the image of the holonomy action lies in the subgroup $\H$.
 
 The purpose of this text is to define rigorously the concepts above and to construct a classifying object and a universal $\H$-quasi-fibration sequence, of the form
 $$F\to B(\ast,\aut_\H F,F)\to B\aut_\H F$$
 which classifies $\H$-fibration sequences. Here $\aut_\H F$ will denote the submonoid of $\aut F$ of self-homotopy equivalences $\varphi$ such that $[\varphi]\in \H \leqslant \E(F)$. Note that when $\H=\E(F)$ is the whole group, then $\aut_\H F=\aut F$ and we recover the usual classifying space and the original result of Stasheff \cite{sta}. The other extreme case $\H=1$, makes $B\aut_{\H}F$ weakly equivalent to the universal cover of $B\aut F$. The rational homotopy type of this last space is well-known (see for instance \cite[\S VII]{tan}). Other choices of the subgroup $\H$ yield intermediate situations between the two described above.
 
  Although this classification may be part of the folklore, we are not aware of any detailed reference in the literature and the purpose of this note is to fill this gap in a rigorous way. The resulting classifying space is a well known and interesting object, see for instance \cite{drorza}. Particularly interesting examples and applications of this classification, in the rational homotopy category, are being considered in \cite{fefuenmu}; there it is shown that for a good choice of the subgroup $\H$, the classifying space has the homotopy type of a rational space. Furthermore, its classifying properties are fundamental for obtaining homotopical information about such space.
  
   As a final remark, I stress that this result differs from the classical classification of topological $G$-bundles for a given topological group $G$ (see for instance \cite{BGM})  in the same extent as the classification of any  class of fibrations differs from the corresponding bundle counterpart.
  
  The text is structured as follows: in the first section, the concepts of holonomy action and $\H$-fibration sequence are introduced. In the second, the geometric-bar construction is used to prove the classification Theorem \ref{main2}. Finally, in the third part, we give a pointed version of these results.

\section*{Acknowledgement}
I would like to thank Scherer J\'er\^ome for his useful advice and corrections.

\section{$\H$-fibration sequences}

Throughout this text, any considered topological space is compactly generated and  weakly Hausdorff. Also, the fiber $F$ and the base $B$ shall be always considered connected spaces of the homotopy type of a CW-complex and $b_0$ is a fixed basepoint in $B$ such that its inclusion is a cofibration. By a {\em fibration sequence}, we mean a sequence
$$
F\stackrel{\omega}{\longrightarrow}E\stackrel{\pi}{\longrightarrow} B
$$
where $p$ is a Hurewicz fibration, and the map $\omega\colon F\stackrel{\simeq_w}{\longrightarrow}F_0=\pi^{-1}(b_0)$ is a weak homotopy equivalence.
 A map between two fibration sequences, $F\stackrel{\omega_1}{\longrightarrow}E_1\stackrel{\pi_1}{\longrightarrow}X$ and $F\stackrel{\omega_2}{\longrightarrow}E_1\stackrel{\pi_2}{\longrightarrow}X$, is a homotopy commutative diagram of the form,
$$
\begin{tikzcd}
& E_1 \arrow{dr}{\pi_1} \arrow{dd}{f} &\\
F\arrow{ur}{\omega_1} \arrow{dr}[swap]{\omega_2} && B\\
& E_2\arrow{ur}[swap]{\pi_2}
\end{tikzcd}
$$
In particular, $f$ is a weak homotopy equivalence and there is no loss of generality if the right triangle is imposed to be strictly commutative. 
 The maps of fibration sequences generate a equivalence relation (by imposing  symmetry and associativity) and we denote by $\fib(B,F)$ the corresponding quotient set.

\begin{remark}
 Note that this definition of $\fib(B,F)$ (which is the one considered in \cite{Allaud}, for example) differs with the one considered in \cite{sta} or \cite{may}, where a fibration with fiber $F$ is a Hurewicz fibration $\pi:E\to B$ whose fiber is weakly equivalent to $F$, but no explicit weak equivalence is given. These two approaches give rise to different equivalence classes and, therefore, to  different classifying spaces. However, both spaces are weakly equivalent (see \cite[explanation below Theorem 1.2]{May80}). Since we are considering base spaces which are of the homotopy type of a CW-complex, both approaches coincide in our setting.
\end{remark}

Recall how the fundamental group of the base acts on the homotopy automorphism of the fiber.
Given a fibration sequence $F\stackrel{\omega}{\longrightarrow}E\stackrel{\pi}{\longrightarrow}B$ and $\beta\in \pi_1(B,b_0)$, the lifting property of the fibration yields an induced homotopy equivalence $\bar\beta\colon F_0=\pi^{-1}(b_0)\stackrel{\simeq}{\to} F_0$. Together with the bijection $\omega_*\colon[F,F]\stackrel{\cong}{\to}[F,F_0]$ produces the homotopy equivalence $\hat\beta={\omega_*}^{-1}[\bar\beta \circ \omega]\in[F,F]$. That is, $\omega\circ \hat\beta \simeq \bar\beta\circ \omega$. Here $[X,Y]$ indicates free homotopy classes of (non-based) maps from $X$ to $Y$.

 Let $\aut F$ be the topological monoid of self homotopy equivalences of $F$ and write ${\cal E}(F)=\pi_0(\aut F)$. Fix a subgroup $\H\leqslant {\cal E}(F)$ and consider $\aut_\H F\subset\aut F$ the submonoid of homotopy equivalences $\varphi$ such that $ [\varphi]\in \H$.
 
 In particular, the process described above, gives, for each fibration sequence $F\stackrel{\omega}{\longrightarrow}E\stackrel{\pi}{\longrightarrow}B$, a map
 $$\pi_1(B,b_0)\to \E(F),\;\;\; [\beta]\mapsto [\hat\beta]$$
 which will be called the \emph{holonomy action of $\pi_1(B,b_0)$ on the fiber}.

\begin{definition}\label{fibh} A fibration sequence $F\stackrel{\omega}{\longrightarrow}E\stackrel{p}{\longrightarrow}B$ is an $\H$-fibration sequence if the image of the holonomy action is contained in $\H$. This means that
 $\hat\beta\in\H$ for any $\beta \in\pi_1(B,b_0)$. 
\end{definition} 
A direct inspection shows the following proposition.
\begin{proposition} Given two fibration sequences over $B$ with fiber $F$ and a map of fibration sequences between them, if one of them is an $\H$-fibration sequence, the other one also is.
\end{proposition}

In particular, we can consider $\Fib_\H(B,F)$
 the set of equivalence classes of $\H$-fibration sequences over $B$ with fiber $F$, which is a subset of $\Fib(B,F)$.

\begin{example}
In the extreme case, $\H=\E(F)$, we recover the usual notion of fibration sequence, since $\Fib_\H(B,F)=\Fib(B,F)$.  In the other extreme case, $\H=1$, an $\H$-fibration sequence is a fibration sequence with a trivial holonomy action. As an intermediate case, we can take the homology functor (or more generally, any functor from the homotopy category of topological spaces),  and take $\H$ as those automorphisms which are the identity in the homology.

\end{example}
 We finish this section recalling two usual constructions to fix the notation. 
 
 \begin{definition}
 Given a fiber sequence  $F\stackrel{\omega}{\longrightarrow}E\stackrel{p}{\longrightarrow}B$ and a based map $f:A\to B$, the \emph{pullback fibration sequence} is defined as
$$F\stackrel{f^*\omega}{\longrightarrow} f^*E=\{(x,a)\in E\times A \mid f(a)=\pi(x)\}\stackrel{f^*\pi}{\longrightarrow} A$$
where $f^*\omega(x)=(\omega(x),a_0)$ and $f^*\pi(x,a)=a$.
\end{definition}

It is immediate to check that the pullback fibration sequence of an $\H$-fibration sequence is an $\H$-fibration sequence.

\begin{remark} If $\{a_0\}\hookrightarrow A$ is a cofibration, an arbitrary map $f:(A,a_0) \to (B,b_0)$ is (free) homotopy equivalent to a based map, so we can always assume that the maps between base spaces are based.
\end{remark}

\begin{definition} For an arbitrary map $\pi:E\to B$, we write
$$\Gamma E=\{(x,\beta)\in E\times B ^I\mid \beta(0)=\pi(x)\}$$
$$\Gamma \pi:\Gamma E\to B,\;\;\; (x,\beta)\mapsto \beta(1)$$
where $B^I$ is the space of paths on $B$. The map $\Gamma \pi$ is a fibration with \emph{homotopy fiber}
$$\F=(\Gamma \pi)^{-1}(b_0)=\{(x,\beta)\in E\times B^I \mid \beta(0)=\pi(x),\; \beta(1)=b_0\}$$

If the natural homotopy equivalence $i:E\to \Gamma E,\; i(x)=(x,c_{\pi(x)})$, where $c_{\pi(x)}$ is the constant path at $\pi(x)$, induces a weak homotopy equivalence $\pi^{-1}(b)\to (\Gamma\pi)^{-1}(b)$ for all $b\in B$, we say that $\pi:E\to B$ is a \emph{quasi-fibration}.

By a \emph{quasi-fibration sequence} we mean a sequence
$$F \stackrel{\omega}{\longrightarrow} E \stackrel{\pi}{\longrightarrow} B$$
such that $\pi$ is a quasi-fibration and $i\circ \omega:F\to \F$ is a weak homotopy equivalence. The \emph{associated fibration sequence} is the fibration sequence
$$F\stackrel{i \circ \omega}{\longrightarrow} \Gamma E\stackrel{\Gamma \pi}{\longrightarrow} B$$

Finally, we say that a quasi-fibration sequence is an \emph{$\H$-quasi-fibration sequence} if its associated fibration sequence is an $\H$-fibration sequence.

\end{definition}

\section{The classification theorem}

We strongly rely in the classical reference \cite{may} from which we recall some facts. Let $G$ be a topological monoid with the identity element $e$ a strongly nondegenerate basepoint, and let $X$ and $Y$ be left and right $G$-spaces respectively. The {\em geometric bar construction  $B(Y,G,X)$} is the geometric realization of the simplicial set whose space of $j$ simplices is $Y\times G^j\times X$ and the face and degeneracy operators are given by,
$$
d_i (y,g_1,\dots, g_j, x)=\left\{ \begin{array}{cl}
(y\cdot g_1,g_2,\dots, g_j,x) & \text{ if } i=0,\\
(y,g_1,\dots,g_{i-1},g_i \cdot g_{i+1},g_{i+2},\dots,g_j,x) & \text{ if } 1\leq i<j,\\
(y,g_1,\dots,g_{j-1},g_j\cdot x) & \text{ if } i=j,\\
\end{array}\right.
$$
$$
s_i(y,g_1,\dots, g_j,x)=(y,g_1,\dots, g_i,e,g_{i+1},\dots, g_j,x).
$$

This is a functor from the category of triples $(Y,G,X)$ as above, whose morphisms are triples $(g,f,h)$ where $f\colon G\to G'$ is a map of monoids and $g\colon Y\to Y'$ and $h:X\to X'$ are $f$-equivariant maps. Of particular interest is the map
$$
p\colon B(*,G,X)\to B(*,G,*)=BG
$$
induced by the trivial $G$-map $X\to *$.
Then if $G$ is a grouplike topological monoid, i.e. $\pi_0(G)$ is a group, then \cite[Theorem 7.6]{may} asserts that $p$ is a quasifibration with fiber $X$. In particular, the inclusion in the homotopy fiber
$$i:X\stackrel{\simeq _w}{\longrightarrow} \F=(\Gamma p)^{-1}(\ast)=\{ (z,\beta)\mid z\in B(\ast,G,X),\beta:[0,1]\to BG, p(z)=\beta(0),\beta(1)=\ast\}$$
is a weak homotopy equivalence, where $i(x)=(x,c_\ast)$, $\ast\in BG$ is the unique point in the unique 0-simplex and $x\in B(\ast,G,X)$ lies in the 0-simplex X.

In particular, whenever
 $F$ is of the homotopy type of a CW-complex, the {\em universal fibration sequence}
$$
F \stackrel{i}{\longrightarrow}
 \Gamma B(\ast,\aut F,F) \stackrel{\Gamma p}{\longrightarrow} B\aut F
$$
classifies fibrations with fiber $F$. Explicitly, see \cite[Theorem 9.2]{may}, for any space $B$ of the homotopy type of a CW-complex, the map,
$$
\Lambda\colon [B,B\aut F]\stackrel{\cong}{\longrightarrow}\fib(B,F), \quad\Lambda [f]=f^*\Gamma p
$$
is a natural bijection.
Given $\H\leqslant\quie(F)$ a  subgroup, consider the map
$$
p_\H\colon B(\ast,\aut_\H F,F)\to B(\ast,\aut_\H F,\ast)=B\aut_\H F.
$$

Finally, let's recall that whenever $G$ is grouplike, then , there is a natural weak homotopy equivalence ( see \cite[Proposition 8.7]{may})
$$\zeta:G\xrightarrow{\simeq_w} \Omega BG,\;\;\; g\mapsto \gamma_g$$
where $\gamma_g:[0,1]\to BG$ is a path such that $\gamma_g(t)=(g,t)$ where this element lies in the 1-skeleton of $BG$.

With these elements we can finally prove the classification theorem.

\begin{theorem}\label{main2}
The map,
$$
\Lambda_\H\colon [B,B\aut_\H F]\stackrel{\cong}{\longrightarrow}\fib_\H(B,F), \quad\Lambda_\H [f]=f^*\Gamma p_{\H},
$$
is a natural bijection which fits in the following commutative diagram
\begin{equation}\label{diagram}
\begin{tikzcd}
 { [B,B\aut F]  } \arrow[]{r}[swap]{\cong} \arrow{r}{\Lambda} & \fib(B,F)  \\
{[B,B\aut_\H F] }\arrow[]{u}{(Bj)_*}\arrow[]{r}[swap]{\cong} \arrow{r}{\Lambda_\H} & \fib_\H(B,F).\arrow[hook]{u}
\end{tikzcd}
\end{equation}
\end{theorem}
Here the left vertical map is induced by the inclusion of monoids $j\colon \aut_\H F\hookrightarrow \aut F$.

\begin{proof}  {\em (i)} For any $[f]\in [X,B\aut_\H F]$, $\Lambda_\H [f]$ is indeed an $\H$-fibration sequence:
For each $g\in \aut_\H F$, lift the path $\gamma_g=\zeta(g)$, in the fibration 
$$
\Gamma p_\H:\Gamma B(\ast,\aut_\H F,F)\to B\aut_\H F$$
 to obtain a homotopy equivalence $\bar\gamma_g:\F\to \F$ of the homotopy fiber. Then a straightforward inspection shows that the following diagram commutes up to homotopy,
$$\begin{tikzcd}
F \arrow{r}{g} \arrow{d}{i}[swap]{\simeq_w}  & F\arrow{d}{i}[swap]{\simeq_w}\\
{\cal F}\arrow{r}{\bar\gamma_g} & {\cal F}.
\end{tikzcd}$$
In particular, we obtain that $[\hat\gamma_g]=[g]\in \E(F)$. Since $\pi_1(B\aut_\H F)\cong\pi_0(\aut_\H F)\cong\H$, we have that the holonomy action factorizes through 
$$\pi_1(B\aut_\H(F))\cong \H\xrightarrow{\id} \H\subset\E(F)$$
so its image lies in $\H$. Thus
 $F\stackrel{i}{\longrightarrow} B(\ast,\aut_\H F,\ast) \xrightarrow{\Gamma p_\H}  B\aut_\H F$ is an $\H$-fibration sequence.  Thus, the pullback $\Lambda_\H [f]$ is an $\H$-fibration sequence.

\medskip

{\em (ii)} The diagram (\ref{diagram}) commutes:

It is enough to show that the fibrations with fiber $F$,
$$\Gamma p_\H\colon \Gamma B(\ast, \aut_\H F,F)\to B\aut_\H F$$
and
$$(Bj)^* \Gamma p\colon (Bj)^*\Gamma B(\ast,\aut F,F)\to B\aut_\H F$$
are equivalent.

For it, define the map 
$$\mu:\Gamma B(\ast,\aut_\H F,F)\to \Gamma B(\ast,\aut F, F),\;\;\; (z,\beta)\mapsto (B(\ast,j,\id_F)(z), Bj\circ\beta)$$

It can be verified that this is a well-defined map. Furthermore,
applying \cite[Proposition 7.8]{may} to the map $j:\aut_\H F\to \aut F$, we get that $\mu$ is map of fibrations, this means that, restricted to each fiber, $\mu$ is a weak homotopy equivalence.

Now consider the following commutative diagram:
$$
\begin{tikzcd}
\Gamma B(\ast,\aut_\H F,F) \arrow[bend left]{drr}{\mu} \arrow[bend right]{ddr}{\Gamma  p_\H}
\arrow[dashed]{dr}{}
\\
& (Bj)^*\Gamma B(\ast,\aut F,F) \arrow{r} \arrow{d}{(Bj)^*\Gamma p}
&\Gamma B(\ast,\aut F,F) \arrow{d}{\Gamma p}\\
&B\aut_\H F \arrow{r}{Bj} & B\aut F.
\end{tikzcd}
$$
where the dashed arrow exists because of the universal property of the pullback. Both $\mu$ and the upper horizontal arrow are maps of fibrations, so also is the dashed arrow. Therefore we have obtained a map of fibrations, between our two initial fibrations.

\medskip

{\em (iii)} The map $\Lambda_\H$ is injective:

Equivalently, we will show that $(Bj)_*:[B,B\aut_\H F]\to [B,B \aut F]$ is injective. By \cite[Remark 8.9]{may} the map $Bj\colon B\aut_\H F\to B\aut F$ is equivalent to a quasifibration with fiber $B(\aut F,\aut_\H F,\ast)$. A long exact sequence argument, using that $\pi_n(Bj)\colon \pi_n(B\aut_\H F)\to \pi_n(B\aut F)$ is an isomorphism for $n\geq 2$, shows that $B(\aut F,\aut_\H F,\ast)$ is weakly equivalent to the discrete space $\quie(F)/\H$.

Using the CW approximation theorem, we can find CW-complexes $Y,Z$ and weak homotopy equivalences $\kappa\colon Y\to B\aut_\H F$ and $\lambda\colon Z\to B\aut F$ such that
 $$\begin{tikzcd}
 Y \arrow{r}{\kappa} \arrow{r}[swap]{\simeq_w} \arrow{d}{\xi} & B\aut_\H F \arrow{d}{Bj}\\
 Z \arrow{r}{\lambda} \arrow{r}[swap]{\simeq_w} & B\aut F\\
 \end{tikzcd}
 $$
 commutes up to homotopy and $\xi$ is a covering map. The associated subgroup of the covering map is given by the image of $\H\leqslant {\cal E}(F)$ under the isomorphisms $\quie(F)\cong \pi_1(B\aut F)$ and $\pi_1(\lambda)\colon \pi_1(Z)\stackrel{\cong}{\to} \pi_1(B\aut F)$. Therefore we have a commutative diagram
\begin{equation}\label{diagram2}
\begin{tikzcd}
{[B,Z]} \arrow{r}{\lambda_*} \arrow{r}[swap]{\cong}&
 { [B,B\aut F]  }  \\
 {[B,Y]} \arrow{r}{\kappa_*} \arrow{r}[swap]{\cong} \arrow[]{u}{\xi_*}&
{[B,B\aut_\H F] }\arrow[]{u}{(Bj)_*}
\end{tikzcd}
\end{equation}
in which, by the lifting property of covering maps, $\xi_*\colon[B,Y]\to [B,Z]$ is injective, so is $(Bj)_*:[X,B\aut_\H F]\to [X,B\aut F]$.

\medskip

{\em (iv)} The  map $\Lambda_\H$ is surjective:

We start with an $\H$-fibration sequence $F\stackrel{\omega}{\longrightarrow}E \stackrel{\pi}{\longrightarrow} B$; via $\Lambda^{-1}$ we get a map $f\colon B\to B\aut F$ which, with the notation in (\ref{diagram2}), produces another map $f'\colon B\to Z$ with $\lambda \circ f'\simeq f$. It is then enough to show that $\im \pi_1(f)\subset \H$. In this case, by the lifting property of  covering maps, there exists $\tilde f\colon B\to Y$ such that $\xi \circ \tilde f=f'$. Therefore, $(Bj)\circ \kappa \circ \tilde f\simeq f$ and $\Lambda_\H[\kappa \circ  \tilde f]$ is the original $\H$-fibration sequence.

To finish, we show that in fact, $\im \pi_1(f)\subset \H$. To give an explicit description of $f$, recall from \cite[p.49]{may} the existence of a commutative diagram of fibrations with fiber $\aut F$
\begin{equation}\label{diagram3}
\begin{tikzcd}
PE\arrow{d}{P\pi} & B(PE,\aut F,\aut F) \arrow{l} \arrow{r}{} \arrow{d}{} & B(\ast, \aut F,\aut F)\arrow{d}{}\\
B \arrow[ shift right=1.0ex, dashed]{r}[swap]{ \varphi}  & B(PE,\aut F,\ast) \arrow{l} \arrow{r}{q} & B\aut F\\
\end{tikzcd}
\end{equation}
where:
\begin{itemize}
\item The space $PE$ is the space of  maps $\psi\colon F\to E$ such that $\psi(F)\subset \pi^{-1}(b)$ for some $b\in B$, and $P\pi$ is defined by $P\pi(\psi)=\pi(\psi(F))$, see \cite[Definition 4.3]{may} for more details.

\item  $B(*,\aut F,\aut F)$  is contractible and  all the maps in the right square are induced by the functor $B(-,-,-)$.

\item All the arrows pointing left are weak homotopy equivalences, so, since $B$ is of the homotopy type of a CW-complex, there exists  $\varphi$, a right homotopy inverse.

\end{itemize}

Then, $\Lambda^{-1}[\pi]=[f]$ with $f=q\circ\varphi$. To check that fundamental group sends this map into $\H$, we first choose
a basepoint $B(PE,\aut F,\aut F)$ which, in turn, determines basepoints in any of the spaces in (\ref{diagram3}).  Among the space $PE\times \aut F$ of $0$-simplices of $B(PE,\aut F,\aut F)$ we fix $(\omega,\id_F)$ with
$\omega\colon F\stackrel{\simeq_w}{\to} F_0=\pi^{-1}(b_0)$ the  weak homotopy equivalence given in the data of the fibration sequence. With this choice, the fiber of $P\pi$ is the space $\FF$ of all weak homotopy equivalences $F\stackrel{\simeq_w}{\to} F_0$.

With this choice of base points the long homotopy exact sequences associated to the $\aut F$-fibrations in (\ref{diagram3}) yield a commutative diagram,
$$\begin{tikzcd}
\vdots \arrow{d} &\vdots \arrow{d} & \vdots \arrow{d}\\
\pi_1(PE)  \arrow{d}{\pi_1(P\pi)} &   \pi_1 B(PE,\aut F,\aut F) \arrow{l} \arrow{d}{} \arrow{r }{} & 0 \arrow{d}\\
\pi_1(B) \arrow{d}{\delta} & \pi_1(B(PE,\aut F,\ast)) \arrow{l}[swap]{\cong} \arrow{d} \arrow{r}{\pi_1(q)} & \pi_1(B\aut F) \arrow{d}{\cong}\\
\pi_0(\FF) \arrow{d} & \quie(F) \arrow{l}[swap]{\omega_*} \arrow{r}{\id} \arrow{d}& \quie(F)\arrow{d}\\
\vdots & \vdots & \vdots
\end{tikzcd}
$$
Hence, $\im\pi_1(f)=\im\pi_1(q\circ \varphi)\in\H$ if and only if $\im (w_*^{-1}\circ\delta)\in\H$. However, an easy inspection shows that for any $\beta\in\pi_1(B)$, $\delta(\beta)=\bar\beta \circ\omega\simeq \omega\circ \hat\beta$. Since we have started with an $\H$-fibration sequence, then $\hat\beta\in \H$, so $\omega_*^{-1}\circ \delta$ sends a path $\beta$ to an element in $\H$, which concludes the proof.
\end{proof}

\begin{remark}\label{remark}
As in the ordinary case, a more convenient expression of the universal $\H$-quasi-fibration sequence can be given.
Fix a base point $x_0\in F_0$ whose inclusion is a cofibration and consider the {\em evaluation fibration}
$$
\aut_\H^* F\hookrightarrow \aut_\H F\stackrel{\ev}{\longrightarrow} F
$$
where $\aut_\H^* F$ is the submonoid of $\aut_\H F$ of
 self homotopy equivalences which fix the base point and
 $\ev(g)=g(x_0)$. This yields a weak homotopy equivalence $F\simeq_w \aut_\H F/\aut_\H^* F=B(\aut_\H,\aut_\H^*,*)$ which, in turn, produces equivalent fibrations with fiber $F$,
 $$\Gamma B(\ast,\aut_\H F,F)\to B\aut_\H F\quad\text{and}\quad \Gamma B(\ast,\aut_\H F,\aut_\H F/\aut_\H^* F)\to B\aut_\H F.
 $$
On the other hand, by
 \cite[Remark 8.9]{may}, the quasifibrations
  $$
  B(\ast,\aut_\H F,\aut_\H F/\aut_\H^*F)\to B\aut_\H F \quad\text{and}\quad B\aut_\H^* F\to B\aut_\H F
  $$
  are equivalent.
That is, the fibration sequences,
$$
F\to B\aut_\H^* F\to B\aut_\H F\quad\text{and}\quad
F\to B(\ast,\aut_\H F,F)\to B\aut_\H F
$$
are equivalent.
\end{remark}

Therefore, Theorem \ref{main2} can be reformulated in the following way:
$$F\to B\aut^*_\H F\to B\aut_\H F$$
is the universal $\H$-quasi-fibration sequence which classifies $\H$-fibration sequences.

\section{Based fibrations}\label{section:2}

The previous results can be adapted to the based case. Let's give a brief exposition of the concepts needed for a rigorous definition.

A fibration sequence $F\stackrel{\omega}{\longrightarrow}E\stackrel{\pi}{\longrightarrow} B$  is {\em based} if $\pi$ has a section $\sigma\colon B\to E$ and the weak equivalence $\omega$ sends the basepoint of $F$ to $\sigma(b_0)$. Consider the `whisker construction' (see \cite[Addenda]{may}),
$$\widetilde \pi  \colon\widetilde E\longrightarrow B\quad\text{where}\quad \widetilde E=\frac{E\sqcup (B\times[0,1])}{\sigma(b)\sim(b,0),\forall b \in X},\quad \widetilde \pi(y)=\pi(y),\quad  \widetilde\pi(b,t)=b.
$$
for $y\in E, b\in F$ and $t\in[0,1]$.
This is a fibration with fiber $\widetilde\pi^{-1}(x_0)=\widetilde F_0=F_0\vee [0,1]$.
Given $\beta\in\pi_1(B,b_0)$ the induced homotopy equivalence $\widetilde\beta\colon \widetilde F_0\stackrel{\simeq}{\to} \widetilde F_0$ can be taken to send 1 to 1. Since the inclusion of the basepoint in $F$ is a cofibration, there is a pointed weak homotopy equivalence  $\tilde\omega:F\to \tilde F_0$, sending $x_0$ to 1, which induces a bijection $\widetilde\omega_*\colon [F,F]^*\stackrel{\cong}{\to} [F,\tilde{F_0}]^*$ between pointed homotopy classes of base-preserving maps. This yields the pointed homotopy equivalence $\dot\beta=\widetilde\omega_*^{-1}[\widetilde\beta f]^*\in[F,F]^*$.

 Let $\aut^* F$ be the topological monoid of pointed self homotopy equivalences of $F$ and denote by  ${\cal E}^*(F)=\pi_0(\aut^*F)$. Fix a subgroup $\H\leqslant {\cal E}^*(F)$ and consider $\aut^*_\H F\subset\aut^*  F$ the submonoid of pointed homotopy equivalences $\varphi$ such that their (pointed) homotopy classes $ [\varphi]^*$ belong to $\H$.

\begin{definition} A based fibration sequence $F\stackrel{\omega}{\longrightarrow}E\stackrel{\pi}{\longrightarrow} B$ is a based $\H$-fibration sequence if $\dot\beta\in\H$ for any $\beta \in\pi_1(B,b_0)$. Denote by $\fib^*_\H(B,F)$ the set of equivalence classes of based $\H$-fibrations sequences over $B$ with fiber $F$.
\end{definition}

Again, pullbacks preserve based $\H$-fibration sequences and this definition is independent of the equivalence class of the given based fibration sequences

 Consider the maps
$$p \colon B(\ast,\aut^* F,F)\to B(\ast,\aut^* F,\ast)=B\aut^*F$$
and
$$
p_\H\colon B(\ast,\aut^*_\H F,F)\to B(\ast,\aut^*_\H F,\ast)=B\aut^*_\H F,
$$
both endowed with  the section induced by the inclusion of the basepoint in $F$, and
let $\Lambda\colon [B,B \aut^*F] \stackrel{\cong}{\to}\fib^*(B,F)$ be the natural bijection given in \cite[Theorem 9.2 (b)]{may}.

Then, by requiring all the maps and fibrations involved in the proof of  Theorem $\ref{main2}$ to be basepoint preserving and based respectively, we obtain:

\begin{theorem}\label{based}
The map,
$$
\Lambda_\H\colon [B,B\aut^*_\H F]\stackrel{\cong}{\longrightarrow}\fib^*_\H(B,F), \quad\Lambda_\H [f]=f^*\widetilde{\Gamma p_{\H}},
$$
is a natural bijection which fits in the following commutative diagram
$$
\begin{tikzcd}
 { [B,B\aut^* F]  } \arrow[]{r}[swap]{\cong} \arrow{r}{\Lambda} & \fib^*(B,F)  \\
{[B,B\aut^*_\H F] }\arrow[]{u}{}\arrow[]{r}[swap]{\cong} \arrow{r}{\Lambda_\H} & \fib^*_\H(B,F).\arrow[hook]{u}
\end{tikzcd}
$$
\end{theorem}

\textsc{
Departamento de \'Algebra, Geometr\'ia y Topolog\'ia, Universidad de M\'alaga, Campus de Teatinos, s/n, 29071, M\'alaga, Spain.}

{\em E-mail address}: \texttt{mario.fuentes@uma.es}

 \end{document}